\documentclass[12pt]{article}
\usepackage{amssymb}
\usepackage{amsmath}
\usepackage[dvips]{graphicx}
\usepackage[T2A]{fontenc}
\usepackage[cp1251]{inputenc}

\usepackage[english,ukrainian]{babel}

\usepackage[cp1251]{inputenc}
\usepackage[ukrainian]{babel}

\usepackage{color}

\Ukrainian

\usepackage{amsthm}

\begin{document}

\noindent{ \bf Volodymyr P. Shchedryk  }

\bigskip

\noindent{ email: shchedrykv@ukr.net }

\bigskip

\noindent {Pidstryhach Institute for Applied Problems of Mechanics and Mathematics, National Academy of Sciences of Ukraine}

\bigskip

  \noindent \textbf{A greatest common divisor and  a least common multiple of solutions of  a linear matrix equation}

\bigskip

\noindent {\it We describe the explicit form of a   left greatest common divisor and  a  least common multiple of solutions of a    solvable linear matrix equation over  a commutative    elementary divisor domain.  We  prove that these left greatest common divisor and   least common multiple are also solutions of the same  equation.}

 \bigskip

{\bf Introduction}.
The investigation  of solutions  of linear matrix   equations was started  in the second half of the XIX century by  Sylvester. Many analytical,  approximate methods and algorithms for finding solutions have been developed. The concepts of the  Moore–Penrose inverse  and the Kronecker  product were actively used.
The question of separating of solutions with certain prescribed  properties  was also considered.
In particular, solutions  with the symmetry condition were studied in [1-3], Hermitian positively defined solutions in [4,5], solutions with a minimal rank in [6], diagonal and triangular solutions in [7].

In the present  article, we investigate solutions of the linear matrix   equation $ BX = A $ over a commutative    elementary divisor domain [8]. Note that a commutative    elementary divisor domain  is a commutative domain over which each matrix is equivalent to a diagonal matrix, each diagonal element (invariant factor) divides the next one.
Examples of   elementary divisor rings are Euclidean rings,    principal ideal rings,  adequate rings, a ring of formal power series over a field of rational numbers with a free integer term (see  [9, 10]).

According to [11, Theorem 2, p. 218] the equation $ BX = A $ has a solution if and only if the corresponding  invariant factors of the matrices $ B $ and $\left[
\begin{array}{cc}
  A & B
\end{array}
\right] $
coincide. In particular, solutions  of a  solvable equation can be found using  the method proposed in  [11, Theorem 1, p. 215].

Let  $ BX = A $ be a solvable equation. Therefore there exists a matrix $ C $ such that $ A = BC. $ It means that the matrix $ B $ is a left divisor of the matrix $ A $. The contrary statement is also correct.   Consequently, the matrix equation $ BX = A $ is solvable if and only if the matrix $ B $ is a  left divisor of the matrix $ A $. Therefore, the solutions of this equation are  left quotients   of   $ A $  by   $ B $. The structure of these quotients are described in [12, 13].
A quotient   of  $ a $ by $ b $ in $ R $ is denoted as $ \frac {a} {b}\in R$.
However, already  the symbol $ \frac {A} {B} $ in matrix rings over $ R $ does not make sense,  because  this quotient is not always  uniquely determined.
Consequently,  the question of choice of that solutions of a   matrix equation $ BX = A $  which are   "generating" others    is arises.
Since  a left greatest common divisor   of matrices over  a commutative    elementary divisor domain (see definition below) defined uniquely  up to  right associativity (see [12, Theorem 1.11]) we can define $ \frac {A} {B} $ as a greatest common divisor of the solutions of a matrix equation $ BX = A $ .
Logically, the question arises: whether a  left g.c.d. of  solutions of    $ BX = A $  is again a solution. In the proposed article we receive an affirmative  answer.
At the same time, we also show that a left least common multiple of solutions  of this equation has the same property.

\bigskip {\bf Main result}.
If $A = BC$, then   $B$ is called a left divisor
 of  $A,$ and
 $A$ is called a right
multiple of  $B.$ If $A=DA_{1} $ and $F=DB_{1}$, then   $D$ is called a left common  divisor of   $A$ and $F.$ In addition, if   $D$ is a right multiple of each left common   divisor of
$A$ and $F$, then   $D$ is called a left  greatest   common  divisor ({\bf g.c.d.})  of   $A$ and $F$.

If $S=MT=NK$, then the matrix $S$ is  called a right common  multiple  of  $M$ and $N$. Moreover, if the matrix $S$ is a right divisor of each right common  multiple of   $M$ and $N$, then $S$ is called a left least common multiple ({\bf l.c.m.})   of  $M$ and $N$.

Our main result is the following.
\bigskip

\noindent
{\bf Theorem}. {\it
Let  $R$ be  a commutative    elementary divisor domain. If \linebreak $BX=A$ is   a solvable matrix equation over $R$, where $A,X,B\in {\rm M}_n(R)$ then a  left g.c.d. and  a left l.c.m. of solutions  of this equation    are also  solutions of $BX=A$. }

\bigskip

\noindent {\bf Proof}. By definition  of the elementary divisor domain,  the matrices $A, B$ can be written as
 $A  =P^{-1} {\rm E} Q^{-1}$ and $B= V^{-1} \Phi U^{-1}$,  in which
$$
{\rm E} = \mathrm {diag}(\varepsilon_1, \ldots , \varepsilon_k, 0,
\ldots , 0),  \qquad \varepsilon_i  | \varepsilon_{i+1}, \qquad i= 1, \ldots , k-1,
$$
$$
\Phi =\mathrm {diag}(\varphi_1, \ldots , \varphi_t, 0, \ldots ,
0), \qquad  \varphi_j  | \varphi_{j+1},\qquad  j= 1, \ldots , t-1,
$$
 and   $P, V, Q, U \in {\rm GL}_n(R).$
The matrix $B$ is a left divisor of  $A$  (see  [Theorem 4.1, 12]), if and only if  $ V=L P,$ where $L $ belongs to  the set
$$
{\bf L}(\text{\rm E}, \Phi)= \{ L \in {\rm GL}_n(R) {\rm |}\;\exists\ S \in M_n(R) :\ L \text{\rm E} = \Phi S \}.
$$
Moreover,   $ t \geq k$ and  $\varphi_i | \varepsilon_i$ $(i= 1 , \ldots , k)$ by  [Theorem 4.6 , 12].
Furthermore,
\begin{equation}\label{BBfg}
{\bf L} (\text{\rm E}, \Phi)=
\left\{L\in {\rm GL}_n(R) \mid L:=\left[
\begin{smallmatrix}
  L_1 & L_3 \\
  L_2 & L_4 \\
  \bf 0 & L_5
\end{smallmatrix}
\right]\right\}
\end{equation}
(see [Theorem 4.5, 11]), in which $L_3, L_4, L_5$ are arbitrary matrices of appropriate sizes,
\[
L_1=\left[
\begin{array}{ccccc}
  l_{11} &  l_{12} & \ldots &  l_{1.k-1} &  l_{1k} \\
   \frac{\varphi_2}{(\varphi_2,\varepsilon_1)}l_{21} & l_{22} & \ldots & l_{2.k-1} & l_{2k} \\
  \ldots & \ldots & \ldots & \ldots & \ldots \\
   \frac{\varphi_k}{(\varphi_k,\varepsilon_1)}l_{k1} &
   \frac{\varphi_k}{(\varphi_k,\varepsilon_2)}l_{k2} & \ldots
   &  \frac{\varphi_k}{(\varphi_k,\varepsilon_{k-1})}l_{k.k-1} & l_{kk}
\end{array}
\right],
\]
\[
L_2=\left[
\begin{array}{ccc}
   \frac{\varphi_{k+1}}{(\varphi_{k+1},\varepsilon_{1})}l_{k+1.1}  & \ldots &
    \frac{\varphi_{k+1}}{(\varphi_{k+1},\varepsilon_{k})}l_{k+1.k} \\
  \ldots & \ldots & \ldots \\
  \frac{\varphi_t}{(\varphi_t,\varepsilon_1)}l_{t1} & \ldots &
  \frac{\varphi_t}{(\varphi_t,\varepsilon_k)}l_{tk}
\end{array}
\right].
\]
Consequently,  $B=(LP)^{-1}\Phi U^{-1}.$
By definition of ${\bf L} (\text{\rm E}, \Phi)$,  we have
\[
L{\rm E}= \Phi S, \quad \text{so}\quad   {\rm E} = L^{-1} \Phi S.
\]
Now
\[
\begin{split}
A  =P^{-1} {\rm E} Q^{-1}&= P^{-1} ( L^{-1} \Phi S) Q^{-1}= (L P)^{-1} \Phi S  Q^{-1}\\
&=((L P)^{-1} \Phi U^{-1})( U S  Q^{-1})=BC,
\end{split}
\]
where $C:=U S  Q^{-1}$.
So a matrix $C$ is the solution of the equation $BX=A.$

If  $ C_1 \not=C$ is    also solution (i.e. $ BC_1 = A$),   then
$$
A=BC_1=BC, \quad \text{so}\quad  B(C_1-C)={\bf 0}.
$$
Hence, $F:=C_1-C \in Ann_r(B)$ and $C_1 \in \{C + Ann_r(B)\}$, where  $ Ann_r (B) $ is the   right annihilator of   $B$. Moreover, if  $C_2 \in \{C + Ann_r(B)\}$, then  $C_2 =C +N$ for some   $N \in Ann_r(B)$, and
$$
BC_2=B(C+N)=BC+BN=A+ {\bf 0}=A.
$$
Consequently, the set $\{C + Ann_r (B)\} $ consists of all solutions of the equation $ BX = A$.

Let us describe   the form of elements  of the set   $\{C + Ann_r(B)\}$.   Using  [Theorem   1.15, 12] and the fact that $B= V^{-1} \Phi U^{-1}$, we obtain that
\[
 Ann_r (B)=\Big\{U\left[ \begin{smallmatrix} {\bf 0}_{t \times n} \\ {D} \end{smallmatrix}\right]\mid D\in M_{(n-t)\times n}(R)\Big\}.
\]
Note that  $C:=U S  Q^{-1}$.  Using  the same idea as in  the proof of  [Theorem 4.5, 12], we obtain that $
S=\left[
\begin{smallmatrix}
  M_1 & {\bf 0} \\
  M_2 & {\bf 0} \\
  M_3 & M_4
\end{smallmatrix}
\right]$,
in which
$$ M_1=\left[
\begin{array}{ccccc}
  \frac{\varepsilon_1}{\varphi_1}l_{11} &  \frac{\varepsilon_2}{\varphi_1}l_{12} &
 \ldots &  \frac{\varepsilon_{k-1}}{\varphi_1}l_{1.k-1} &  \frac{\varepsilon_k}{\varphi_1}l_{1k} \\
   \frac{\varepsilon_1}{(\varphi_2,\varepsilon_1)}l_{21} & \frac{\varepsilon_2}{\varphi_2}l_{22}
   & \ldots & \frac{\varepsilon_{k-1}}{\varphi_2}l_{2.k-1} & \frac{\varepsilon_{k}}{\varphi_2}l_{2k} \\
  \ldots & \ldots & \ldots & \ldots & \ldots \\
   \frac{\varepsilon_1}{(\varphi_k,\varepsilon_1)}l_{k1} & \ldots & \ldots
   &  \frac{\varepsilon_{k-1}}{(\varphi_k,\varepsilon_{k-1})}l_{k.k-1} & \frac{\varepsilon_{k}}{\varphi_k}l_{kk}
\end{array}
 \right],
 $$
 $$
M_2=\left[
\begin{array}{ccc}
   \frac{\varepsilon_{1}}{(\varphi_{k+1},\varepsilon_{1})}l_{k+1.1}  & \ldots &
    \frac{\varepsilon_{k}}{(\varphi_{k+1},\varepsilon_{k})}l_{k+1.k} \\
  \ldots & \ldots & \ldots \\
  \frac{\varepsilon_1}{(\varphi_t,\varepsilon_1)}l_{t1} & \ldots &
  \frac{\varepsilon_k}{(\varphi_t,\varepsilon_k)}l_{tk}
\end{array}
\right],
 $$
 and
$\left[ \begin{smallmatrix} M_3 & M_4 \end{smallmatrix}\right]$ is an arbitrary matrix from  $  {\rm M}_{(n-t)\times n}(R)$.
Consequently,
\begin{equation}\label{Dffrl}
\begin{split}
C+Ann_r(B)&=USQ^{-1}+U \textstyle\left[ \begin{smallmatrix} {\bf 0}_{t \times n} \\ {D} \end{smallmatrix}\right]\\
&=U\left(S+ \textstyle\left[ \begin{smallmatrix} {\bf 0}_{t \times n} \\ {D}Q \end{smallmatrix}\right] \right)Q^{-1}
=
U\left(S+ \left[ \begin{smallmatrix} {\bf 0}_{t \times n} \\ D_1 \end{smallmatrix}\right] \right)Q^{-1}.
\end{split}
\end{equation}
Using the fact that  $Q$ is an invertible matrix, it is easy to check that
\[
{\rm M}_{(n-t)\times n}(R)Q= {\rm M}_{(n-t)\times n}(R),
\]
so  $D_1:=DQ$  (see  \eqref{Dffrl})  can be  any  matrix from ${\rm M}_{(n-t)\times n}(R)$.
Consequently,
$$
C+Ann_r(B)= U\left( \left[
\begin{smallmatrix}
  M_1 & {\bf 0} \\
  M_2 & {\bf 0} \\
  M_3 & M_4
\end{smallmatrix}
\right] +
\left[ \begin{smallmatrix} {\bf 0}_{t \times n} \\ D_1 \end{smallmatrix}\right] \right)Q^{-1}
= U \left[
\begin{smallmatrix}
  M_1 & {\bf 0} \\
  M_2 & {\bf 0} \\
  T_3 & T_4
\end{smallmatrix}
\right] Q^{-1},
$$
where $\left[ \begin{smallmatrix} T_3 & T_4 \end{smallmatrix}\right]$ is an arbitrary matrix from  ${\rm M}_{(n-t)\times n}(R)$.
It follows that
\begin{equation}\label{VVffdGG}
F:=
  U\left[
\begin{smallmatrix}
  M_1 & {\bf 0} \\
  M_2 & {\bf 0} \\
  {\bf 0} & I_{n-t}
\end{smallmatrix}
\right] Q^{-1} \in \{C+Ann_r(B)\},
\end{equation}
where $I_{n-t}$ is the identity matrix of order $n-t$.
Since $\{C+Ann_r(B)\}$ consists of all solutions of $BX=A$, and
\[
  U\left[
\begin{smallmatrix}
  M_1 & {\bf 0} \\
  M_2 & {\bf 0} \\
T_3 & T_4
\end{smallmatrix}
\right] Q^{-1}=
\left(  U\left[
\begin{smallmatrix}
  M_1 & {\bf 0} \\
  M_2 & {\bf 0} \\
  {\bf 0} & I_{n-t}
\end{smallmatrix} \right]
Q^{-1}\right)
\left(  Q\left[
\begin{smallmatrix}
I_t & {\bf 0} \\
T_3 & T_4
\end{smallmatrix} \right]
Q^{-1}\right)=FM,
\]
the solution  $F$ of the equation $BX=A$ is a left divisor of all elements from $\{C+Ann_r(B)\}$.
Using the fact that $ F \in \{C + Ann_r (B)\} $, we obtain that  $F$ is a left g.c.d. of $\{C + Ann_r (B)\} $.

Consider the matrix
\begin{equation}\label{SujjJg}
N:=  U\left[
\begin{smallmatrix}
  M_1 & {\bf 0} \\
  M_2 & {\bf 0} \\
{\bf 0} & {\bf 0}
\end{smallmatrix}
\right] Q^{-1}=
\left(  U\left[
\begin{smallmatrix}
I_t & {\bf 0} \\
{\bf 0} & {\bf 0}
\end{smallmatrix} \right]
U^{-1}\right)
\left(  U\left[
\begin{smallmatrix}
  M_1 & {\bf 0} \\
  M_2 & {\bf 0} \\
T_3 & T_4
\end{smallmatrix} \right]
Q^{-1}\right).
\end{equation}
Obviously,   $N\in \{C + Ann_r (B)\}$ and it is a left multiple of all elements of $\{C + Ann_r (B)\}$, so $N$ is a  left l.c.m.  of $\{C + Ann_r (B)\}$. The proof is done.
{\hfill \rm  $\Box$}

\bigskip

\noindent {\bf Corollary 1}. {\it
Each left g.c.d. and each  left l.c.m. of solutions of a solvable  equation $BX = A$ have the following forms:
$$
 U\left[
\begin{smallmatrix}
  M_1 & {\bf 0} \\
  M_2 & {\bf 0} \\
{\bf 0} & I_{n-t}
\end{smallmatrix}
\right] Q^{-1}, \qquad
 U\left[
\begin{smallmatrix}
  M_1 & {\bf 0} \\
  M_2 & {\bf 0} \\
{\bf 0} & {\bf 0}
\end{smallmatrix}
\right] Q^{-1},
$$
respectively.} {\hfill \rm  $\Box$}

\bigskip

\noindent {\bf Corollary 2}. {\it The set
\[
\left\{
\left(  U\left[
\begin{smallmatrix}
  M_1 & {\bf 0} \\
  M_2 & {\bf 0} \\
  {\bf 0} & I_{n-t}
\end{smallmatrix} \right]
Q^{-1}\right)
\left(  Q\left[
\begin{smallmatrix}
I_t & {\bf 0} \\
T_3 & T_4
\end{smallmatrix} \right]
Q^{-1}\right)\quad \mid \quad \left[ \begin{smallmatrix}  T_3 & T_4 \end{smallmatrix}\right]\in {\rm M}_{(n-t)\times n}(R) \right\}
\]
consists of all solutions of the equation $ BX = A. $}
{\hfill \rm  $\Box$}

\bigskip

\noindent {\bf Example}. Let $$
A:=\mathrm {diag} (1,\; 2, \; 6,\; 0,\; 0,\; 0, \;0), \;\;
 B:=
\left[
\begin{matrix}
  0 & 1 & 0 & 0 & 0 & 0 & 0  \\
    0 & -2 & 2 & 0 & 0 & 0 & 0  \\
      1 & 0 & 0 & 0 & 0 & 0 & 0  \\
        0 & 0 & 0 & 0 & 0 & 0 & 0  \\
          -2 & 0 & -12 & 0 & 12 & 0 & 0  \\
            0 & 0 & 0 & 0 & 0 & 0 & 0  \\
              -2 & 0 & -4 & 4 & 0 & 0 & 0
\end{matrix}
\right]
$$
are matrices over  $\mathbb{Z}$.
It is easy to check that the Smith forms of matrices $A$ and $B$ are $E:=A$ and $\Phi:=\mathrm {diag} (1,\; 1, \; 2,\; 4,\; 12,\; 0, \;0)$, respectively.  Moreover, $B=V^{-1} \Phi$  in which
$$ V^{-1}=
\left[
\begin{matrix}
  0 & 1 & 0 & 0 & 0 & 0 & 0  \\
    0 & -2 & 1 & 0 & 0 & 0 & 0  \\
      1 & 0 & 0 & 0 & 0 & 0 & 0  \\
        0 & 0 & 0 & 0 & 0 & 0 & 0  \\
          -2 & 0 & -6 & 0 & 1  & 0 & 0  \\
            0 & 0 & 0 & 0 & 0 & 0 & 0  \\
              -2 & 0 & -2 & 1 & 0 & 0 & 0
\end{matrix}
\right].
$$
Consider the matrix  equation $BX=A$ over $\mathbb{Z}$.
Using \eqref{BBfg} and the fact that  $\Phi | {\rm E}$,  we get that  ${\bf L}(\text{\rm E}, \Phi)\neq \emptyset $  by [Theorem 4.5, 11] and
$$
{\bf L}(\text{\rm E}, \Phi)=
\left\{\left[
\begin{array}{ccc|cccc}
  * & * & * & * & * & * & *  \\
  * & * & * & * & * & * & *  \\
  2h_{31} &* & * & * & * & * & *  \\ \hline
  4h_{41}& 2h_{42} & 2h_{43} & *  & *  & *  & *   \\
 12h_{51} & 6h_{52} & 2h_{53} &  *  & *  & *  & *  \\ \hline
   0 & 0 & 0 & *  & *  & *  & *  \\
   0 & 0 & 0 & *  & *  & *  & *
\end{array}
\right]\right\}.
$$
The matrix $A$ has the form $A=I_7{\rm E}I_7. $ It follows that $ V = LI_7 = L $ (see notation of the Theorem) and
\[ V=L=
\left[
\begin{array}{ccc|cccc}
  0 & 0 & 1 & 0 & 0 & 0 & 0  \\
    1 & 0 & 0 & 0 & 0 & 0 & 0  \\
      2 & 1 & 0 & 0 & 0 & 0 & 0  \\ \hline
        4& 2 & 2 & 0 & 0 & 0 & 1  \\
          12 & 6 & 2 & 0 & 1  & 0 & 0  \\ \hline
            0 & 0 & 0 & 0 & 0 & 1 & 0  \\
              0 & 0 & 0 & 1 & 0 & 0 & 0
\end{array}
\right]\in {\bf L}(\text{\rm E}, \Phi).
\]
This yields that the equation $ BX = A $ has a solution (see  [Theorem 4.1, 12]). Furthermore,  $ L{\rm E}= \Phi S$ for a matrix
\[
S:=
\left[
\begin{array}{ccc|cccc}
  0 & 0 & 6 & 0 & 0 & 0 & 0  \\
    1 & 0 & 0 & 0 & 0 & 0 & 0  \\
      1 & 1 & 0 & 0 & 0 & 0 & 0  \\ \hline
        1& 1 & 3 & 0 & 0 & 0 & 0  \\
          1 & 1 & 1 & 0 & 0  & 0 & 0  \\ \hline
            * & * & * & * & * & * & *  \\
             * & * & * & * & * & * & *
\end{array}\right].
\]
Consequently, a  left g.c.d.  and a  l.c.m. of the solutions of $ BX = A $ are:
\[
\left[
\begin{array}{ccc|cccc}
  0 & 0 & 6 & 0 & 0 & 0 & 0  \\
    1 & 0 & 0 & 0 & 0 & 0 & 0  \\
      1 & 1 & 0 & 0 & 0 & 0 & 0  \\ \hline
        1& 1 & 3 & 0 & 0 & 0 & 0  \\
          1 & 1 & 1 & 0 & 0  & 0 & 0  \\ \hline
            0 & 0& 0 & 0 & 0  & 1 & 0  \\
             0 & 0 & 0 & 0 & 0 & 0 & 1
\end{array}
\right]\quad \text{and}\quad
\left[
\begin{array}{ccc|cccc}
  0 & 0 & 6 & 0 & 0 & 0 & 0  \\
    1 & 0 & 0 & 0 & 0 & 0 & 0  \\
      1 & 1 & 0 & 0 & 0 & 0 & 0  \\ \hline
        1& 1 & 3 & 0 & 0 & 0 & 0  \\
          1 & 1 & 1 & 0 & 0  & 0 & 0  \\ \hline
            0 & 0& 0 & 0 & 0  & 0 & 0  \\
             0 & 0 & 0 & 0 & 0 & 0 & 0
\end{array}
\right],
\]
respectively (see \eqref{VVffdGG} and \eqref{SujjJg}).
Moreover,  these matrices are solutions of the equation $BX=A$.
{\hfill \rm  $\diamond $}

\bigskip
\noindent
{\bf Problem}. Describe  that rings over them a left g.c.d. and  a left l.c.m. of solutions of a solvable matrix equation $BX=A$ are also   solutions.

\begin{enumerate}

\item  \textit{ Vetter W. J.}   Vector structures and solutions of linear matrix equations // Linear Algebra Appl. -- 1975. -- \textbf{10}. --   № 2. --  P. 181-188.

\item  \textit{ Magnus  J. R.,   Neudecker H. }  The elimination matrix: Some lemmas and applications // SIAM. J. on Algebraic and Discrete Methods. 1980. -- \textbf{1}. --   № 4. -- P. 422–449.

\item  \textit{   Don F.J. } On the symmetric solutions of a linear matrix equation // Linear Algebra Appl. -- 1987. -- \textbf{93}. --   P. 1-7.

\item  \textit{   Khatri C.G.,    Mitra S.K.}  Hermitian and Nonnegative Definite Solutions of Linear Matrix Equations // Journal on Applied Mathematics. --  1976. -- \textbf{31}. --  № 4. -- P. 579-585.

\item  \textit{  Ran A.,    Reurings M.}   A Fixed Point Theorem in Partially Ordered Sets and Some Applications to Matrix Equations // Proceedings of the American Mathematical Society. -- 2004.  -- \textbf{132}. --   № 5. -- P. 1435–1443.

\item  \textit{   Recht B.,   Fazel M.,    Parrilo P.}   Guaranteed minimum-rank solutions of linear matrix
equations via nuclear norm minimization // \newline  doi.org/10.1137/070697835.

\item  \textit{   Magnus J.R. } L-structured matrices and linear matrix equations // Linear Algebra Appl. -- 1983.    -- \textbf{14}. --   № 1. --   P. 67-88.

\item  \textit{Kaplansky I}. Elementary divisor and modules // Trans. Amer. Math. Soc. -- 1949. -- \textbf{66}. -- P. 464--491.

\item   \textit{  Bovdi V.A., Shchedryk V.P.} Commutative Bezout domains of stable range 1.5 // Linear Algebra and Appl. -- 2019. -- \textbf{568}. -- P. 127-134.

\item  \textit{  Zabavsky B.V.} Diagonal reduction of matrices over rings. -- Lviv. --
Math Studies Monograph Series. -- VNTL Publishers. --  2012. -- 251 p.

 \item  \textit{  Kazimirsky P.S.} Decomposition of matrix polynomials into factors. Naukova dumka, 1981. -- 224 p.

\item  \textit{ Shchedryk  V.P. } Factorization of matrices over elementary divisor rings. -- Lviv: Pidstryhach Institute for Applied Problems of Mechanics and Mathematics of the NAS of Ukraine, 2017. -- 304 p. \linebreak
    http://www.iapmm.lviv.ua/14/index.htm

\item  \textit{  Shchedryk V.P.}  Factorization of
 matrices over elementary divisor domain   // Algebra  Discrete Math. -- 2009. -- № 2. -- P. 79--99.

\end{enumerate}

\end{document}